\documentclass[a4paper,10pt]{article}

\usepackage{amssymb, amscd, verbatim, amsthm}
\usepackage[all]{xy}
\usepackage{hyperref}
\usepackage{amsmath}

\begin{document}

\title{Abelian varieties over finitely generated fields and the conjecture of Geyer and Jarden on torsion}

\author{Sara Arias-de-Reyna\\
        Institut f{\" u}r Experimentelle\\
        Mathematik,\\
        45326 Essen, Germany\\
        \texttt{\small sara.arias-de-reyna@uni-due.de}
        \vspace{3mm}
\and
        Wojciech Gajda\\
        Department of Mathematics,\\
        Adam Mickiewicz University,\\
        61614 Pozna\'{n}, Poland\\
        \texttt{\small gajda@amu.edu.pl}\\
        \vspace{3mm}
\and
        Sebastian Petersen\footnote{the corresponding author}\\
        Universit\"at Kassel,\\
        Fachbereich Mathematik,\\
        34132 Kassel,
        Germany\\
        \texttt{\small sebastian.petersen@unibw.de}
}
\parindent0em
\parskip1em

\maketitle

\newtheorem{leer}{}[section]
\newtheorem{thm}[leer]{Theorem}
\newtheorem{conj}[leer]{Conjecture}
\newtheorem{defi}[leer]{Definition}
\newtheorem{rema}[leer]{Remark}
\newtheorem{prop}[leer]{Proposition}
\newtheorem{lemm}[leer]{Lemma}
\newtheorem{coro}[leer]{Corollary}
\newtheorem{quest}[leer]{Question}
\newtheorem{claim}[leer]{Claim}

\newcounter{abcSatz}
\newtheorem*{abcSatz}{Main Theorem}
\renewcommand{\theabcSatz}{\Alph{abcSatz}}
\newtheorem*{abcVermutung}{Conjecture of Geyer and Jarden}
\newtheorem*{abcFolgerung}{Corollary}

\newcommand{\prm}{\ {\mathrm{prime}}}

\newcommand{\ilim}{\mathop{\varinjlim}\limits}
\newcommand{\plim}{\mathop{\varprojlim}\limits}

\newcommand{\OO}{{\cal O}}
\newcommand{\PP}{{\cal P}}
\renewcommand{\AA}{{\cal A}}
\newcommand{\BB}{{\cal B}}
\newcommand{\CC}{{\cal C}}
\newcommand{\DD}{{\cal D}}
\newcommand{\LL}{{\cal L}}
\newcommand{\UU}{{\cal U}}
\newcommand{\TT}{{\cal T}}
\newcommand{\MM}{{\cal M}}
\newcommand{\FF}{{\cal F}}
\newcommand{\NN}{{\cal N}}
\newcommand{\RR}{{\cal R}}
\newcommand{\KK}{{\cal K}}
\newcommand{\XX}{{\cal X}}
\newcommand{\YY}{{\cal Y}}
\renewcommand{\SS}{{\cal S}}
\newcommand{\Mon}{{\cal M}}

\newcommand{\Occ}{{\mathrm{Occ}}}
\newcommand{\sep}{_{\mathrm{sep}}}
\newcommand{\bad}{_{\mathrm{bad}}}
\newcommand{\tor}{_{\mathrm{tor}}}
\newcommand{\ab}{_{\mathrm{ab}}}
\newcommand{\Spec}{\mathrm{Spec}}
\newcommand{\Eig}{\mathrm{Eig}}
\newcommand{\dro}{\mathrm{drop}}
\newcommand{\Sym}{\mathrm{Sym}}
\newcommand{\Sp}{\mathrm{Sp}}
\newcommand{\GL}{\mathrm{GL}}
\newcommand{\amp}{\mathrm{amp}}
\newcommand{\Emb}{\mathrm{Emb}}
\newcommand{\GSp}{\mathrm{GSp}}
\newcommand{\PSp}{\mathrm{PSp}}
\newcommand{\SL}{\mathrm{SL}}
\newcommand{\Hom}{\mathrm{Hom}}
\newcommand{\Alg}{\mathrm{Alg}}
\newcommand{\ggT}{{\mathrm{ggT}}}
\newcommand{\can}{{\mathrm{can}}}
\newcommand{\trdeg}{\mathrm{trdeg}}
\newcommand{\Mor}{\mathrm{Mor}}
\newcommand{\mf}{\mathfrak}
\newcommand{\ol}[1]{\overline{#1}}
\newcommand{\Var}{\mathrm{Var}}
\newcommand{\AbVar}{\mathrm{AbVar}}
\newcommand{\Res}{\mathrm{Res}}
\newcommand{\Grp}{\mathrm{Grp}}
\newcommand{\Sch}{\mathrm{Sch}}
\newcommand{\Inert}{\mathrm{Inert}}
\newcommand{\Iso}{\mathrm{Iso}}
\newcommand{\Cur}{\mathrm{Cur}}
\newcommand{\Alb}{\mathrm{Alb}}
\newcommand{\Rat}{\mathrm{Rat}}
\newcommand{\Aut}{\mathrm{Aut}}
\newcommand{\End}{\mathrm{End}}
\newcommand{\DivCor}{\mathrm{DivCor}}
\newcommand{\degr}{\mathrm{deg}}
\newcommand{\ran}{\mathrm{rg}}
\newcommand{\rk}{\mathrm{rk}}
\newcommand{\ke}{\mathrm{ker}}
\newcommand{\coke}{\mathrm{coker}}
\newcommand{\im}{\mathrm{im}}
\newcommand{\chara}{\mathrm{char}}
\newcommand{\rot}{\mathrm{rot}}
\newcommand{\ord}{{\mathrm{ord}}}
\newcommand{\Rest}{\ \mathrm{Rest}\ }
\newcommand\F{\mathbb{F}}
\newcommand\Gal{\mathrm{Gal}}
\newcommand\pr{\mathrm{pr}}
\newcommand\Id{\mathrm{Id}}
\newcommand\Stab{\mathrm{Stab}}
\newcommand\calB{\mathcal{B}}

\newcommand{\dR}[4]{#1\div #2=#3\ \Rest #4}

\newcommand{\TS}{TS}

\def\Pp{{\mathbb{P}}}
\def\Ff{{\mathbb{F}}}
\def\Rr{{\mathbb{R}}}
\def\Cc{{\mathbb{C}}}
\def\Qq{{\mathbb{Q}}}
\def\Zz{{\mathbb{Z}}}
\def\Aa{{\mathbb{A}}}
\def\Nn{{\mathbb{N}}}
\def\Gg{{\mathbb{G}}}

\begin{abstract}
In this paper we prove the Geyer-Jarden conjecture on the torsion part of the Mordell-Weil group for a
large class of abelian varieties defined over finitely generated fields of arbitrary characteristic.
The class consists of all abelian varieties with {\em big monodromy}, i.e., such that the image of Galois representation on $\ell$-torsion points, for almost all primes $\ell,$ contains the full symplectic group.
\end{abstract}

\newpage

\section*{Introduction}
\addcontentsline{toc}{section}{Introduction}
Let $A$ be a polarized abelian variety defined over a finitely generated field $K.$ Denote by $\widetilde{K}$ (respectively, $K\sep$) the algebraic (resp., separable) closure of $K$. It is well known that the Mordell-Weil
group $A(K)$ is a finitely generated $\Zz$-module. On the other hand $A(\widetilde{K})$ is a
divisible group with an infinite torsion part $A(\widetilde{K})\tor$ and $A(\widetilde{K})$ has infinite rank, unless $K$ is algebraic over a finite field.
Hence, it is of fundamental interest to study the structure of the groups $A(\Omega)$ for infinite algebraic
extensions $\Omega/K$ smaller than $\widetilde{K}$. For example, Ribet in \cite{ribet1981} and Zarhin in \cite{zarhin1987}
considered the question of finiteness of $A(K\ab)\tor,$
where $K\ab$ is the maximal abelian extension of $K$.
\medskip

We denote by $G_K:=G(K\sep/K)$ the absolute Galois group of $K$. For a
positive integer $e$ and for
$\sigma=(\sigma_1,\sigma_2,\dots, \sigma_e)$ in the group $G_K^e=$ $G_K{\times} G_K{\times} \dots {\times} G_K,$ we denote by $K\sep(\sigma)$ the subfield in $K\sep$ fixed by $\sigma_1,\sigma_2,\dots, \sigma_e.$
There exists a substantial literature on arithmetic properties of the fields $K\sep(\sigma).$ In particular,
the Mordell-Weil groups $A(K\sep(\sigma))$ have been already studied, e.g., Larsen formulated a conjecture
in \cite{larsen2003} on the rank of $A(K\sep(\sigma))$ (cf. \cite{larsenim2008}, \cite{geyerjarden2006rank} for results supporting the conjecture of Larsen).
\medskip

In this paper we consider the torsion part of the groups $A(K\sep(\sigma)).$
In order to recall the conjecture which is mentioned in the title, we agree to say that a property
${\cal A}(\sigma)$ {\it holds for almost all $\sigma\in G_K^e$}, if ${\cal A}(\sigma)$ holds
for all $\sigma\in G_K^e,$ except for a set of measure zero with respect to the (unique) normalized Haar measure
on the compact group $G_K^e.$ In \cite{geyerjarden1978} Geyer and Jarden proposed the following
conjecture on the torsion part of $A(K\sep(\sigma))$.

\begin{abcVermutung}\label{gjconjecture}
Let $K$ be a finitely generated field. Let $A$ be an abelian
variety defined over $K$.
\begin{enumerate}
\item[a)] For almost all $\sigma\in G_K$ there are
infinitely many prime numbers $\ell$ such that the group $A(K\sep(\sigma))[\ell]$ of $\ell-$division points is nonzero.
\item[b)] Let $e\ge 2$. For almost all $\sigma\in G_K^e$ there are only
finitely many prime numbers $\ell$ such that the group $A(K\sep(\sigma))[\ell]$ of $\ell-$division points is nonzero.
\end{enumerate}
\end{abcVermutung}

It is known due to the work of Jacobson and Jarden \cite{jacobsonjarden2001}
that for all $e\ge 1$, almost all $\sigma\in G_K^e$ and all
primes $\ell$ the group $A(K\sep(\sigma))[\ell^\infty]$ is finite. This was formerly
part (c) of the conjecture. Moreover the Conjecture  is known for
elliptic curves \cite{geyerjarden1978}. Part (b) holds true provided $ \chara(K)=0$ (see
\cite{jacobsonjarden2001}). In a very recent preprint Zywina proves part (a) in the special case
where $K$ is a number field (cf. \cite{zywinapreprint}), stengthening results of
Geyer and Jarden \cite{gj2005}.

As for today, for an abelian variety $A$ of dimension $\ge 2$ defined over a finitely generated field of positive
characteristic, parts (a) and (b) of the Conjecture of Geyer and Jarden are open and part (a) is open over a finitely generated transcendental extension of~$\Qq$.

In this paper we prove the Conjecture of Geyer and Jarden for abelian varieties with big monodromy.
To formulate our main result we need some notation.
Let $\ell\neq\chara(K)$ be a prime number. We denote by
$\rho_{A[\ell]}:G_K\longrightarrow \Aut(A[\ell])$ the Galois representation attached
to the action of $G_K$ on the $\ell$-torsion points of $A.$
We define $\Mon_K(A[\ell]):=\rho_{A[\ell]}(G_K)$ and call
this group {\em the mod-{$\ell$} monodromy group of $A/K$}.
We fix a polarization and denote by $e_\ell{:}A[\ell]\times A[\ell]\to \mu_\ell$ the corresponding
Weil pairing. Then $\Mon_{K}(A[\ell])$ is a subgroup of the group of symplectic similitudes $\GSp(A[\ell],e_{\ell})$ of the Weil pairing.
We will say that $A/K$ {\em has big monodromy} if there exists a constant $\ell_0$ such that
$\Mon_{K}(A[\ell])$ contains the symplectic group $\Sp(A[\ell], e_\ell),$ for every prime number $\ell\geq \ell_0.$
Note that the property of having big monodromy does not depend on the choice of the polarization.

The main result of our paper is the following

\begin{abcSatz}\label{ThmA}[cf. Thm. \ref{GJmain}, Thm. \ref{GJmain2}]
Let $K$ be  a finitely generated field and $A/K$ an abelian variety with big monodromy.
Then the Conjecture of Geyer and Jarden holds true for $A/K$.
\end{abcSatz}

Surprisingly enough, the most difficult to prove is the case (a) of
the Conjecture for abelian varieties with big
monodromy, when $\chara(K)>0$. The method of our proof relies in
this case on the Borel-Cantelli Lemma of measure theory and on a
delicate counting argument in the group $\Sp_{2g}(\Ff_\ell)$ which
was modeled after a construction of subsets $S^{\prime}(\ell)$ in
${\mathrm{SL}}_2(\Ff_\ell)$ in Section 3 of the classical paper
\cite{geyerjarden1978} of Geyer and Jarden.

It is interesting to combine the main Theorem with existing
computations of monodromy groups for certain families of abelian
varieties. Certainly, the most prominent result of this type is the
classical theorem of Serre (cf. \cite{serre1984}, \cite{serre1985}
for the number field case; the generalization to finitely generated
fields of characteristic zero is well-known): {\em If $A$ is an
abelian variety over a finitely generated field $K$ of
characteristic zero with $\End(A)=\Zz$ and $\dim(A)=2, 6$ or odd,
then $A/K$ has big monodromy.} Here
$\End(A)=\End_{\tilde{K}}(A_{\tilde{K}})$ stands for the absolute
endomorphism ring of $A$.

In this paper we focus our attention at abelian
varieties with $\End(A)=\Zz,$ which have been recently considered
by Chris Hall in his open image theorem \cite{hall2008}.
To simplify notation, we will say that an
abelian variety $A$ over a finitely generated field $K$ {\em is of
Hall type}, if $\End(A)=\Zz$ and $K$ has a discrete valuation at
which $A$ has semistable reduction of toric dimension one. The
following result, proven in our paper \cite{agp2011}, gives examples of
abelian varieties with big monodromy in all dimensions (and including the case
${\mathrm{char}}(K)>0$): {\em If $A$ is an abelian variety of Hall
type over a finitely generated infinite field $K$, then $A/K$ has
big monodromy.} In the special case where $K$ is a global field this
has recently been shown by Hall (cf. \cite{hall2007},
\cite{hall2008}). The generalization to an arbitrary finitely
generated ground field $K$ is carried out in \cite{agp2011} using
methods of group theory, finiteness properties of the fundamental group of schemes
and Galois theory of large field extensions. In
combination with the main Theorem we obtain the following
\newpage

\begin{abcFolgerung}\label{ThmD}  
Let $A$ be an abelian variety over a finitely generated infinite
field $K$. Assume that either condition i) or ii) is satisfied.
\begin{enumerate}
\item[i)] $A$ is of Hall type.
\item[ii)] $\chara(K)=0$, $\End(A)=\Zz$ and $\dim(A)=2, 6$ or odd.
\end{enumerate}
Then the Conjecture of Geyer and Jarden holds true for $A/K$.
\end{abcFolgerung}

We thus obtain over every finitely generated infinite field and for every dimension families of abelian
varieties for which the Conjecture of Geyer and Jarden holds true. In the case when
$\chara(K)>0$ the Corollary offers the first evidence for the conjecture of Geyer and Jarden on torsion going beyond the case of elliptic curves.

We warmly thank Gerhard Frey, Dieter Geyer, Cornelius Greither and Moshe Jarden for conversations and useful comments on the topic of this paper.
\medskip

\section{Notation and background material}
In this section we fix notation and gather some background material
on Galois representations that is
important for the rest of this paper.

If $K$ is a field, then we denote by $K\sep$ (resp. $\widetilde{K}$) the
separable (resp. algebraic) closure of $K$ and by $G_K=G(K\sep/K)$ its absolute
Galois group. If $G$ is a profinite (hence compact) group, then it has
a unique normalized Haar measure $\mu_G$. The expression ``assertion ${\cal A}(\sigma)$ holds
for almost all $\sigma \in G$'' means ``assertion ${\cal A}(\sigma)$ holds true for
all $\sigma$ outside a zero set with respect to $\mu_G$''.
A finitely generated field is by definition a field which is
finitely generated over its prime field.
Let $X$ be a scheme of finite type over a field $K$. For a geometric point
$P\in X(\tilde{K})$ we denote by $K(P)\subset \tilde{K}$ the residue field at $P$.

For $n\in\Nn$ coprime to $\chara(K)$,
we let $A[n]$ be the group of $n$-torsion
points in $A(\tilde{K})$ and define $A[n^\infty]=\bigcup\limits_{i=1}^\infty A[n^i]$.
For a prime $\ell\neq \chara(K)$ we denote by $T_\ell(A)=\plim_{i\in\Nn} A[\ell^i]$ the $\ell$-adic
Tate module of $A$. Then $A[n]$, $A[n^\infty]$ and $T_\ell(A)$ are $G_K$-modules in a natural way.

If $M$ is a $G_K$-module (for example $M=\mu_n$ or $M=A[n]$ where
$A/K$ is an abelian variety), then we shall denote
the corresponding representation
of the Galois group $G_K$ by
$$\rho_M: G_K\to \Aut(M)$$
and define $\Mon_K(M):=\rho_M(G_K)$. We define $K(M):=K\sep^{\ker(\rho_M)}$ to be the fixed field in
$K\sep$ of the kernel of $\rho_M$. Then $K(M)/K$ is a Galois extension and
$G(K(M)/K)\cong \Mon_K(M)$.

If $R$ is a commutative ring with $1$ (usually $R=\Ff_\ell$ or $R=\Zz_\ell$) and $M$ is a finitely generated
free $R$-module equipped with a non-degenerate
alternating bilinear pairing $e:M\times M\to R'$ into
a free $R'$-module of rank $1$ (which is a multiplicatively written $R$-module
in our setting
below), then we denote by
$$\Sp(M, e)=\{f\in \Aut_R(M)\quad|\quad\forall x, y\in M: e(f(x), f(y))=e(x, y)\}$$
the corresponding symplectic group and by
$$\GSp(M, e)=\{f\in \Aut_R(M)\quad|\quad \exists\varepsilon\in R^\times:\forall x, y\in M:
e(f(x), f(y))=\varepsilon e(x, y)\}$$
the corresponding group of symplectic similitudes.

Let $n$ be an integer coprime to $ \chara(K)$ and $\ell$ be a prime different from
$ \chara(K)$. Let $A/K$ be an abelian variety. We denote by $A^\vee$ the dual abelian
variety and let $e_n: A[n]\times A^\vee[n]\to \mu_n$ and $e_{\ell^\infty}: T_\ell A\times
T_\ell A^\vee\to \Zz_\ell(1)$ be the corresponding Weil pairings. If $\lambda: A\to A^\vee$
is a polarization, then we deduce Weil pairings
$e_n^\lambda: A[n]\times A[n]\to \mu_n$ and $e_{\ell^\infty}^\lambda:
T_\ell A\times T_\ell A\to \Zz_\ell(1)$ in the obvious way.
If $\ell$ does not divide $\deg(\lambda)$ and if $n$ is coprime to $\deg(\lambda)$,
then $e_n^\lambda$ and $e_{\ell^\infty}^\lambda$ are non-degenerate, alternating,
$G_K$-equivariant pairings. Hence we have representations
$$\rho_{A[n]}: G_K\to \GSp(A[n], e_n^\lambda),$$
$$\rho_{T_\ell A}: G_K\to \GSp(T_\ell A, e_{\ell^\infty}^\lambda)$$
with images
$\Mon_K(A[n])\subset \GSp(A[n], e_n^\lambda)$ and
$\Mon_K(T_\ell A)\subset \GSp(T_\ell A, e_{\ell^\infty}^\lambda)$.
We shall say that an abelian variety $(A, \lambda)$ over a field $K$
has {\sl big monodromy}, if there is a constant $\ell_0>\max( \chara(K), \deg(\lambda))$ such
that $\Mon_K(A[\ell])\supset \Sp(A[\ell], e_\ell^\lambda)$ for every prime number
$\ell\ge \ell_0$.

\section{Properties of abelian varieties with big monodromy}

Let $(A, \lambda)$ be a polarized
abelian variety with big monodromy over a finitely generated field $K$.
Then $\Sp(A[\ell], e_\ell^\lambda)\subset
\Mon_K(A[\ell])$ for sufficiently large primes $\ell$. In this section we 
determine $\Mon_K(A[n])$ completely for every ``sufficiently large'' integer $n$. The
main result (cf. Proposition \ref{fullresult} below) is due to Serre in the number field case, and
the general case requires only a slight adaption of Serre's line of reasoning. However, as the 
final outcome is somewhat different in positive characteristic, we do include the details.
Proposition \ref{fullresult} will be crucial for our results on the Conjecture of Geyer and Jarden.

Now let $K$ be an arbitrary field and $A/K$ an abelian variety.
Recall that for every algebraic extension $L/K$ we defined
$\Mon_L(A[n])=\rho_{A[n]}(G_L)$ ($n$ coprime to $\chara(K)$) and
$\Mon_L(T_\ell A)=\rho_{T_\ell A}(G_L)$ ($\ell>\chara(K)$ a prime number).
Furthermore
the representations induce isomorphisms $G(L(A[n])/L)\cong \Mon_L(A[n])$ and
$G(L(A[\ell^\infty]/L)\cong \Mon_L(T_\ell A)$. Note that
$\Mon_L(T_\ell A)\to \Mon_L(A[\ell^i])$ is surjective (because $G(L(A[\ell^\infty])/L)\to
G(L(A[\ell^i])/L)$ is surjective) for every integer $i$. Clearly $\Mon_L(A[n])$ is
a subgroup of $\Mon_K(A[n])$.

\begin{rema} \label{monrem} If $L/K$ is
a {\em Galois} extension, then $\Mon_L(A[n])$ is a {\em normal} subgroup of $\Mon_K(A[n])$ and
the quotient group $\Mon_K(A[n])/\Mon_L(A[n])$ is isomorphic to $G(L\cap K(A[n])/K)$.
\end{rema}

\begin{prop} \label{monpropnew}Let $K$ be a field and $(A, \lambda)$ a polarized
abelian variety over $K$ with big monodromy. Let $L/K$ be
an abelian Galois extension with $L\supset \mu_\infty$.
Then there is a constant $\ell_0>\max(\chara(K), \deg(\lambda))$ with the following
properties.
\begin{enumerate}
\item[a)] $\Mon_L(T_\ell A)=\Sp(T_\ell A, e_{\ell^\infty}^\lambda)$ for all primes $\ell\ge \ell_0$.
\item[b)] Let $c$ be the product of all prime numbers $\le \ell_0$.
Then $\Mon_L(A[n])=\Sp(A[n], e_n^\lambda)$
for every integer
$n$ which is coprime to $c$.
\end{enumerate}
\end{prop}

{\em Proof.}
Part a). There is a constant $\ell_0 >\max(\chara(K), \deg(\lambda), 5)$ such that
$\Mon_K(A[\ell])\supset \Sp(A[\ell], e_{\ell}^\lambda)$ for all primes $\ell\ge \ell_0$, because $A$ has
big monodromy. Let $\ell\ge \ell_0$ be a prime and
define $K_\ell:=K(\mu_\ell)$. Then basic properties of the Weil pairing
imply that $G(K_\ell(A[\ell])/K_\ell)\cong\Mon_{K_\ell}(A[\ell])=\Sp(A[\ell], e_{\ell}^\lambda)$.
This group is perfect, because $\ell\ge 5$ (cf. \cite[Theorem 8.7]{taylor}).
As $L/K_\ell$ is an abelian Galois extension, $\Mon_L(A[\ell])$ is a normal subgroup of the perfect
group $\Mon_{K_\ell}(A[\ell])$ and the quotient $\Mon_{K_\ell}(A[\ell])/\Mon_L(A[\ell])$ is isomorphic
to a subquotient of $G(L/K)$ (cf. Remark \ref{monrem}), hence abelian. This implies that
$$\Mon_L(A[\ell])=\Mon_{K_\ell}(A[\ell])=\Sp(A[\ell], e_\ell^\lambda).$$
Denote by $p: \Sp(T_\ell A, e_{\ell^\infty}^\lambda)\to \Sp(A[\ell], e_\ell^\lambda)$ the
canonical projection. Then $\Mon_L(T_\ell A)$ is a closed subgroup of $\Sp(T_\ell A, e_{\ell^\infty}^\lambda)$
with $$p(\Mon_L(T_\ell A))=\Mon_L(A[\ell])=\Sp(A[\ell], e_\ell^\lambda).$$ Hence
$\Mon_L(T_\ell A)=\Sp(T_\ell A, e_{\ell^\infty}^\lambda)$ by \cite[Proposition 2.6]{larsen1995}.

Part b). Consider the map
$$\rho: G_L\to \prod_{\ell\ge \ell_0} \Mon_L(T_\ell A)=\prod_{\ell\ge \ell_0} \Sp(T_\ell A, e_{\ell^\infty}^\lambda)$$
induced by the representations $\rho_{T_\ell A}$ and denote by $X:=\rho(G_L)$ its image.
Then $X$ is a closed subgroup of $\prod_{\ell\ge \ell_0} \Sp(T_\ell A, e^\lambda_{\ell^\infty})$.
If $\pr_\ell$ denotes the $\ell$-th projection of the product, then
$\pr_\ell(X)=\Sp(T_\ell A, e^\lambda_{\ell^\infty})$. Hence \cite[Section 7, Lemme 2]{serretovigneras}
implies that
$X=\prod_{\ell\ge \ell_0} \Sp(T_\ell A, e^\lambda_{\ell^\infty})$, i.e. that $\rho$ is {\em surjective}.

Let $c$ be the product of all prime numbers $\le \ell_0$. Let $n$ be an integer coprime to $c$.
Then $n=\prod_{\ell|n\prm} \ell^{v_\ell}$ for certain integers $v_\ell\ge 1$. The canonical map
$r: \Mon_L(A[n])\to\prod_{\ell|n\prm} \Mon_L(A[\ell^{v_\ell}])$ is injective. Consider the diagram
$$\begin{xy}
  \xymatrix{
 G_L\ar[d]\ar[r]^{\rho'} & \prod_{\ell|n} \Mon_L(T_\ell A)\ar[d]\ar@{=}[r] &   \prod_{\ell|n}\Sp(T_\ell A, e_{\ell^\infty}^\lambda)\ar[d] \\
 \Mon_L(A[n]) \ar@^{(->}[r]^r & \prod_{\ell|n} \Mon_L(A[\ell^{v_\ell}])\ar@^{(->}[r] &   \prod_{\ell|n} \Sp(A[\ell^{v_\ell}], e_{\ell^{v_\ell}}^\lambda).
}
\end{xy}$$

The vertical maps are surjective. The horizontal map $\rho'$ is surjective as well, because $\rho$ is surjective.
This implies, that the lower horizontal map
$$\Mon_L(A[n])\to\prod_{\ell|n} \Sp(A[\ell^{v_\ell}], e_{\ell^{v_\ell}}^\lambda)$$
is in fact bijective. It follows from the Chinese Remainder Theorem that
the canonical map
$$\prod_{\ell|n} \Sp(A[\ell^{v_\ell}], e_{\ell^{v_\ell}}^\lambda)\to \Sp(A[n], e_n^\lambda)$$
is bijective as well. Assertion b) follows from that.\hfill $\Box$

\begin{coro} \label{extendbigmonodromy}Let $K$ be a field and $(A, \lambda)$ a polarized
abelian variety over $K$ with big monodromy. Then there is a constant $c$ coprime to
$\deg(\lambda)$ and to $\chara(K)$, if $\chara(K)$ is positive, with the following
property: $\Mon_K(A[n])\supset \Sp(A[n], e_n^\lambda)$ for every integer $n$ coprime to $c$.
\end{coro}

{\em Proof.} Let $L=K\ab$ be the maximal abelian extension.
Then there is a constant $c$ as above, such that $\Mon_L(A[n])= \Sp(A[n], e_n^\lambda)$
for every $n$ coprime to $c$ by Proposition \ref{monpropnew}. Furthermore $\Mon_L(A[n])\subset \Mon_K(A[n])$
by the discussion before Remark \ref{monrem}.\hfill $\Box$

Let $K$ be a field and $(A, \lambda)$ a polarized
abelian variety over $K$ with big monodromy.
There is a constant $c$ (divisible by $\deg(\lambda)$ and by $ \chara(K)$, if $\chara(K)\neq 0$) such that
$$\Sp(A[n], e^\lambda_{n})\subset \Mon_K(A[n])\subset
 \GSp(A[n], e^\lambda_{n})$$
 for all $n\in \Nn$ coprime to $c$ (cf. Corollary \ref{extendbigmonodromy}). One can easily determine $\Mon_K(A[n])$
 completely, if $K$ is finitely generated.
 Let $K_{n}:=K(A[n])$. There is a commutative
 diagram
 $$\begin{xy}
  \xymatrix{
 0\ar[r] &  G(K_{n}/K(\mu_{n})) \ar[r]\ar[d] &  G(K_{n}/K)\ar[r]\ar[d]_{\rho_{A[n]}} &G(K(\mu_{n})/K) \ar[r]\ar[d]_{\rho_{\mu_n}} & 0 \\
 0\ar[r] & \Sp(A[n], e^\lambda_{n})\ar[r] & \GSp(A[n], e^\lambda_{n})\ar[r]^\varepsilon & (\Zz/n)^\times\ar[r] & 0
}
\end{xy}$$
with exact rows and injective vertical maps, where $\rho_{\mu_n}$ is the cyclotomic character and $\varepsilon$ is the multiplicator map.
The left hand vertical map is an isomorphism for every $n\in \Nn$ coprime to $c$.
Hence
$$\Mon_K(A[n])=\{f\in \GSp(A[n], e^\lambda_{n})\ |\ \varepsilon(f)\in \im(\rho_{\mu_n})\}.$$

Assume from now on that $K$ is finitely generated.
Then the image of the cyclotomic character involved above has a well known explicit description. Denote by $F$ the
algebraic closure of the prime field of $K$ in $K$ and define $q:=q(K):=|F|\in\Nn\cup \{\infty\}$.
Then, after possibly replacing $c$ by a larger constant, we have
$$
\im(\rho_{\mu_n})=\left\{\begin{array}{ll}
\langle \ol{q} \rangle & \mbox{char}(K)\neq 0,\\
(\Zz/n)^\times & \mbox{char}(K)= 0.\end{array}\right.$$
for all $n\in\Nn$ coprime to $c$. Here $\langle \ol{q} \rangle$
is the subgroup of $(\Zz/n)^\times$
generated by the residue class $\ol{q}$ of $q$ modulo $n$,
provided $q$ is finite. If $q$ is finite,
then we define
$$\GSp^{(q)}(A[n], e^\lambda_{n})=\{f\in \GSp(A[n], e^\lambda_{n})\ |\
\varepsilon(f)\in \langle \ol{q} \rangle\}.$$
Finally we put $\GSp^{(\infty)}(A[n], e^\lambda_{n})=\GSp(A[n], e^\lambda_{n})$.
We have shown:

\begin{prop} \label{fullresult} Let $K$ be a finitely generated field and $(A, \lambda)$
a polarized abelian variety over $K$ with big monodromy.
Let $q=q(K)$. Then there is a constant $c$ (divisible by $\deg(\lambda)$ and by $ \chara(K)$, if $\chara(K)\neq 0$) such that
$\Mon_K(A[n])=\GSp^{(q)}(A[n], e^\lambda_{n})$ for all $n\in \Nn$ coprime
to $c$.
\end{prop}

\section{Proof of the Conjecture of Geyer and Jarden, part b)}

Let $(A, \lambda)$ be a polarized abelian variety of dimension $g$ over a field $K$. In this section
we will use the notation $K_\ell:=K(A[\ell])$ and $G_\ell:=G(K_\ell/K)$ for every prime $\ell\neq  \chara(K)$.
Our main result in this section is the following theorem.

\begin{thm}\label{GJmain} If $(A, \lambda)$ has big monodromy, then for all $e\ge 2$ and almost
all $\sigma\in G_K^e$ (in the sense of the Haar measure)
there are only finitely many primes $\ell$ such that
$A(K\sep(\sigma))[\ell]\neq 0$.
\end{thm}

The following Lemma \ref{MVlemm} is due to
Oskar Villareal (private communication). We thank him for his kind permission to include it into our manuscript.

\begin{lemm} Assume that $A$ has big monodromy. Then there is a
constant $\ell_0$ such that $[K(P):K]^{-1}\le [K_\ell:K]^{-\frac{1}{2g}}$
for all primes $\ell\ge \ell_0$ and all $P\in A[\ell]$, where $K(P)$ denotes the
residue field of the point $P.$\label{MVlemm}
\end{lemm}

{\em Proof.}
By assumption on $A,$ there is a constant $\ell_0$ such that $\Sp(A[\ell], e_\ell^\lambda)\subset
\Mon_K(A[\ell])$ for all primes $\ell\ge \ell_0$. Let $\ell\ge \ell_0$ be a prime and $P\in A[\ell]$.
Then the $\Ff_\ell$-vector space generated inside $A[\ell]$
by the orbit $X:=\{f(P): f\in\Mon_K(A[\ell])\}$
is the whole of $A[\ell]$, because $A[\ell]$ is a simple
$\Ff_\ell[\Sp(A[\ell], e_\ell^\lambda)]$-module. Thus we can choose
an $\Ff_\ell$-basis $(P_1,\cdots, P_{2g})$ of $A[\ell]$ with
$P_1=P$ in such a way that each $P_i\in X$. Then each $P_i$ is
conjugate to $P$ under the action of $G_K$ and $[K(P):K]=[K(P_i):K]$ for
all $i$. The field $K_\ell$ is the composite field $K_\ell=K(P_1)\cdots K(P_{2g})$.
It follows that $$[K_\ell:K]\le [K(P_1):K]\cdots [K(P_{2g}):K]=[K(P):K]^{2g}.$$
The desired inequality follows from that.\hfill $\Box$

The following notation will be used in the sequel: For sequences
$(x_n)_n$ and $(y_n)_n$ of
positive real numbers we shall write $x_n\sim y_n$, provided the
sequence $(\frac{x_n}{y_n})$ converges to a positive real number.
If $x_n\sim y_n$ and $\sum x_n<\infty$, then $\sum y_n<\infty$.

The proof of Theorem \ref{GJmain} will make heavy use of the
following classical fact.

\begin{lemm} {\bf (Borel-Cantelli, \cite[18.3.5]{friedjarden})}\label{borelcantelli}
Let $(A_1, A_2, \cdots)$ be a sequence of measurable subsets of
a profinite group $G$. Let
$$A:=\bigcap_{n=1}^\infty \bigcup_{i=n}^\infty A_i=\{x\in G: \mbox{$x$ belongs to infinitely many $A_i$}\}.$$
\begin{enumerate}
\item[a)] If $\sum_{i=1}^\infty \mu_G(A_i)<\infty$, then $\mu_G(A)=0$.
\item[b)] If $\sum_{i=1}^\infty \mu_G(A_i)=\infty$ and $(A_i)_{i\in\Nn}$ is a $\mu_G$-independent sequence (i.e.  for every finite set $I\subset \Nn$ we have
$\mu_G(\bigcap_{i\in I} A_i)=\prod_{i\in I} \mu_G(A_i)$),
then $\mu_G(A)=1$.
\end{enumerate}
\end{lemm}

{\em Proof of Theorem \ref{GJmain}.} Assume that $A/K$ has big
monodromy and let $\ell_0$ be a constant as in the definition of the
term ``big monodromy''. We may assume that $\ell_0\ge \mbox{char}(K)$.
Let $e\ge 2$ and define
$$X_\ell:=\{\sigma\in G_K^e: A(K\sep(\sigma))[\ell]\neq 0\}$$
for every prime $\ell$. Let $\mu$ be the normalized Haar measure
on $G_K^e$.
Theorem \ref{GJmain} follows from Claim 1 below, because Claim 1
together with the Borel-Cantelli Lemma \ref{borelcantelli} implies that
$$\bigcap_{n\in\Nn}\bigcup_{\ell\ge n\ {\mathrm{prime}}} X_\ell$$
has measure zero.

{\bf Claim 1.} The series $\sum_{\ell\ {\mathrm{prime}}} \mu(X_\ell)$
converges.

Let $\ell\ge \ell_0$ be a prime number.
Note that
$$X_\ell=\bigcup_{P\in A[\ell]\setminus\{0\}} \{\sigma\in G_K^e\ |\ \sigma_i(P)=P\ \mbox{
for all $i$}\}=
\bigcup_{P\in A[\ell]\setminus\{0\}} G_{K(P)}^e.$$
Let $\Pp(A[\ell])=(A[\ell]\setminus\{0\})/\Ff_\ell^\times$ be the projective space
of lines in the $\Ff_\ell$-vector space $A[\ell]$. It is a projective space of
dimension $2g-1$. For $P\in A[\ell]\setminus\{0\}$ we denote by $\ol{P}:=\Ff_\ell^\times P$ the
equivalence class of $P$ in $\Pp(A[\ell])$. For $\ol{P}\in\Pp(A[\ell])$  and $P_1, P_2\in \ol{P}$
there is an $a\in \Ff_\ell^\times$ such that $P_1=aP_2$ and $P_2=a^{-1} P_1$, and
this implies $K(P_1)=K(P_2)$.
It follows that we can write
$$X_\ell=\bigcup_{\ol{P}\in \Pp(A[\ell])} G_{K(P)}^e.$$

Hence
$$\mu(X_\ell)\le\sum_{\ol{P}\in \Pp(A[\ell])} \mu(G_{K(P)}^e)=
\sum_{\ol{P}\in \Pp(A[\ell])} [K(P):K]^{-e},$$
and Lemma \ref{MVlemm} implies
$$\mu(X_\ell)\le\sum_{\ol{P}\in \Pp(A[\ell])} [K_\ell:K]^{-e/2g}=\frac{\ell^{2g}-1}{\ell-1}
[K_\ell:K]^{-e/2g}=\frac{\ell^{2g}-1}{\ell-1}
|G_\ell|^{-e/2g}.$$
But $G_\ell$ contains $\Sp_{2g}(\Ff_\ell)$ and
$$s_\ell:=|\Sp_{2g}(\Ff_\ell)|=\ell^{g^2}\prod_{i=1}^g (\ell^{2i}-1)$$
(cf. \cite{taylor}).
It is thus enough to prove the following

{\bf Claim 2.} The series $\sum_{\ell\ge \ell_0\ {\mathrm{prime}}}
\frac{\ell^{2g}-1}{\ell-1} s_\ell^{-e/2g}$ converges.

But $s_\ell\sim \ell^{g^2+2+4+\cdots+2g}=\ell^{2g^2+g}$ and
$\frac{\ell^{2g}-1}{\ell-1}\sim \ell^{2g-1}$, hence
$$\frac{\ell^{2g}-1}{\ell-1} s_\ell^{-e/2g}\sim \ell^{2g-1} \ell^{-e(g+\frac{1}{2})}=
\ell^{(2-e)g-(1+\frac{e}{2})}\le \ell^{-2},$$
because $e\ge 2$. Claim 2 follows from that. \hfill $\Box$

\section{Special sets of symplectic matrices}
This section contains a construction of certain special sets of symplectic matrices (cf. Theorem \ref{Proposicion} below) that
will play a crucial role in the proof of part a) of the Conjecture of Geyer and Jarden.

Let $g\geq 2$, and let $V$  be a vector space of dimension $2g$ over
a prime finite field $\F_{\ell}$, endowed with a symplectic form
$e: V\times V\rightarrow \F_{\ell}$. Fix a symplectic basis
$E=\{e_1, \dots, e_{2g}\}$ of $V$ such that the symplectic form is
given by the matrix
$$J_g=\begin{pmatrix}J_1 & \ & \ & \ \\
                      \ & J_1 & \ & \ \\
                      \ & \ & \ddots & \ \\
                      \ & \ & \ & J_1\end{pmatrix} \text{ where }J_1=\begin{pmatrix}0 & 1 \\ -1 & 0 \\\end{pmatrix}.$$

For each $A\in \GSp_{2g}(\F_{\ell})$ there is an element $\lambda\in
\F_{\ell}^\times$ such that $e(Av, Aw)=\lambda e(v,  w)$ for all $v,
w\in V$. We will say that the value $\lambda=\varepsilon(A)$
of the multiplicator map $\varepsilon$ is the
\emph{multiplier} of $A$, and we will denote by
$\GSp_{2g}(\F_{\ell})[\lambda]$  the set of matrices in
$\GSp_{2g}(\F_{\ell})$ with multiplier $\lambda$.

\begin{rema} Here we collect some notation. Let $p$ be a prime, $q$ a power of $p$ and $n\in \Nn$.
\begin{itemize}
\item For $n$ not divisible by $p$, we will denote by $\ord_nq$ the order
    of $q$ modulo $n$.
\item Denote by $\GSp^{(q)}_{2g}(\Zz/n\Zz)$ the set of matrices
    in $\GSp_{2g}(\Zz/n\Zz)$ with multiplier equal to a power of $q$ modulo $n$.
\item Denote by $\GSp^{(\infty)}_{2g}(\Zz/n\Zz):=\GSp_{2g}(\Zz/n\Zz)$.
\item Let $\alpha_3, \alpha_4, \dots, \alpha_{2g}, \beta\in
    \F_{\ell}$. Call $u_{\alpha}=e_2 + \alpha_3 e_3 + \cdots+
    \alpha_{2g}e_{2g}$. We denote by $T_{u_{\alpha}}[\beta]$
    the morphism $v\mapsto v + \beta e (v, u_{\alpha}) u_{\alpha}$ (which is a transvection if $\beta\not=0$).
\end{itemize}
\end{rema}

We begin with two easy lemmas that will be essential for Definition
\ref{Definicion}.

\begin{lemm}\label{FixE1} Let $\ell$ be a prime number.
For each $\lambda\in \F_{\ell}^\times$, the matrices of
$\GSp_{2g}(\F_{\ell})[\lambda]$ that fix the vector $e_1$ are of the form
\begin{equation}\label{A}
\left(\begin{array}{c |c | c c c} 1 & d & b_1 & b_2  & \dots\\
\hline 0 & \lambda & 0 & 0 & \dots\\ \hline
0 & d_1 & \ & \ & \ \\
\vdots & \vdots & \ & B & \ \\
\vdots & \vdots & \ & \ & \ \\ \end{array}\right)\end{equation}

with $B=(b_{ij})_{i, j=1, \dots, 2g-2}\in
\GSp_{2g-2}(\F_{\ell})[\lambda]$, $d, d_1, \dots, d_{2g-2}\in
\F_{\ell}$ and
\begin{equation}\label{b_k}b_k=\frac{1}{\lambda}\left( \sum_{j=1}^{g-1} (d_{2j-1} b_{2j, k}- d_{2j}b_{2j-1,
k})\right) \text{ for each }k=1, \dots, 2g-2.\end{equation}

\end{lemm}

{\em Proof.} Let $A\in \GSp_{2g}(\F_{\ell})[\lambda]$ be such that $Ae_1=e_1$.
Let us write the matrix of $A$ with respect to the symplectic basis
$\{e_1, e_2, \dots, e_{2g-1}, e_{2g}\}$. Since $e(e_1, e_k)=0$ for
all $k=3, \dots, 2g$, we obtain that $e(e_1, Ae_k)=0$. Therefore we
can write the matrix $A$ as
\begin{equation*}\left(\begin{array}{c |c | c c c} 1 & d & b_1 & b_2  & \dots\\
\hline 0 & d' & 0 & 0 & \dots\\ \hline
0 & d_1 & \ & \ & \ \\
\vdots & \vdots & \ & B & \ \\
\vdots & \vdots & \ & \ & \ \\ \end{array}\right)\end{equation*}
where in the second row we get all entries zero save the $(2,
2)$-th. Moreover, since $e(e_1, e_2)=1$, we get that $e(e_1,
Ae_2)=e(A e_1, Ae_2)=\lambda e(e_1, e_2)=\lambda$, that is to say,
$d'=\lambda$.

Furthermore, we have that $e(e_2, e_k)=0$ for all $k=3, \dots, 2g$,
hence $e(Ae_2, Ae_k)=0$. These conditions give rise to the
equations \eqref{b_k}. The rest of the conditions one has to impose
imply that $B\in \GSp_{2g-2}(\F_{\ell})[\lambda]$. This proves that
the conditions in the lemma are necessary. On the other hand, one
can check that the product
$$A^t J_g A=\lambda J_g,$$ so they are also sufficient.\hfill $\Box$

\begin{lemm}\label{Independent} The set of matrices in
$\GSp_{2g}(\F_{\ell})[\lambda]$ that do not have the eigenvalue $1$
has cardinality greater than $\beta(\ell, g)\vert \Sp_{2g-2}(\F_{\ell})\vert$, where
$$\beta(\ell, g)=\ell^{2g-1}(\ell^{2g}-1)\frac{\ell-2}{\ell-1}.$$
\end{lemm}

{\em Proof.} The set of  matrices $A\in\GSp_{2g}(\F_{\ell})[\lambda]$ that fix
the vector $e_1$ consists of matrices of the form \eqref{A}, where $B$ belongs
to $\GSp_{2g-2}(\F_{\ell})[\lambda]$ and $b_1, \dots, b_{2g-2}$ are
given by the formula \eqref{b_k} of Lemma \ref{FixE1}. Therefore the cardinality of the set
of such matrices is exactly
$$\ell^{2g-1}\vert \GSp_{2g-2}(\F_{\ell})[\lambda]\vert=\ell^{2g-1}\vert \Sp_{2g-2}(\F_{\ell})\vert.$$
On the other hand, the symplectic group acts transitively on the set
of cyclic subgroups of $V$ (cf. \cite[Thm. 9.9, Ch. 2]{Huppert}). Therefore if a matrix fixes any nonzero
vector, it can be conjugated to one of the above. Hence, to obtain
an upper bound for the number of matrices with eigenvalue 1 one has
to multiply the previous number by the number of cyclic groups of
$V$, namely $\frac{\ell^{2g}-1}{\ell-1}$.

Therefore the set of matrices in $\GSp_{2g}(\F_{\ell})[\lambda]$
that have the eigenvalue $1$ has cardinality less than
$\ell^{2g-1}\frac{\ell^{2g}-1}{\ell-1}\vert \Sp_{2g-2}(\F_{\ell})\vert$. Hence the number of
matrices in $\GSp_{2g}(\F_{\ell})[\lambda]$ that do not have the
eigenvalue $1$ is greater than $\vert \Sp_{2g}(\F_{\ell})\vert -
\ell^{2g-1}\frac{\ell^{2g}-1}{\ell-1}\vert \Sp_{2g-2}(\F_{\ell})\vert$.

Now apply the well known identity (see for instance the proof of
\cite[Theorem 9.3. b)]{Huppert})
\begin{equation}\label{Induccion}\vert \Sp_{2g}(\F_{\ell})\vert =(\ell^{2g}-1)\ell^{2g-1} \vert
\Sp_{2g-2}(\F_{\ell})\vert.\end{equation}

We thus see that the set of matrices in
$\GSp_{2g}(\F_{\ell})[\lambda]$ that do not have the eigenvalue $1$
has cardinality greater than
$\beta(\ell, g)\vert \Sp_{2g-2}(\F_{\ell})\vert$.
\hfill $\Box$

\begin{defi}\label{Definicion} For each $\lambda\in \F_{\ell}^\times$ choose once and for all a subset
$\calB_{\lambda}$ of matrices $B\in \GSp_{2g-2}(\F_{\ell})[\lambda]$
which do not have the eigenvalue 1, with $$\vert
\calB_{\lambda}\vert =\beta(\ell, g-1)\vert
\Sp_{2g-4}(\F_{\ell})\vert$$ (which can be done by Lemma
\ref{Independent}). Define
$$\begin{aligned}&\begin{aligned}S_{\lambda}(\ell)_0:=\{& A \text{ of the shape } \eqref{A} \text{in Lemma \ref{FixE1}  such that:} \\
                 & B\in \calB_{\lambda} \\
                 & d_1, \dots, d_{2g-2}\in  \F_{\ell} \\
                 & d \in \F_{\ell} \setminus \{-(b_1, \dots,
                 b_{2g-2})(\Id-B)^{-1}\begin{pmatrix}d_1, &
\dots, & d_{2g-2}\end{pmatrix}^t\},\\ \end{aligned}\\
& S_{\lambda}(\ell):=\{T_{u_{\alpha}}[\beta]^{-1} \cdot A\cdot
T_{u_{\alpha}}[\beta]: \alpha_3, \dots, \alpha_{2g}, \beta\in
\F_{\ell}, A\in S_i(\ell)_0\},\\\end{aligned}$$

Let $q$ be a power of a prime $p\not=\ell$. Define $$S^{(q)}(\ell):=
\bigcup_{i=1}^{\ord_{\ell}q} S_{q^i}(\ell).$$ Define also
$$S^{(\infty)}(\ell)=\bigcup_{\lambda\in \F_{\ell}^{\times}}
S_{\lambda}(\ell).$$
\end{defi}

\begin{rema} \label{eveins}Clearly $S^{(q)}(\ell)\not=\emptyset$ and $S^{(\infty)}(\ell)\not=\emptyset$. Note
moreover that all matrices in $S^{(q)}(\ell)$ and
$S^{(\infty)}(\ell)$ fix an element of $V$.\end{rema}

This section is devoted to prove the following result.

\begin{thm}\label{Proposicion} The following properties hold:

\begin{itemize}

\item[(1)] Let $q$ be a power of a prime number or $q = \infty$. Then
 $$\sum_{\ell}\frac{\vert
    S^{(q)}(\ell)\vert}{\vert \GSp_{2g}^{(q)}(\F_{\ell})\vert}=\infty.$$
In the first case $\ell$ runs through all prime numbers coprime to
$q$ and in the second case through all prime numbers.

\item[(2)] Let $q$ be a power of a prime number $p$ or $q=\infty$.
Let $\ell_1, \dots, \ell_r$ be different prime numbers. If $q\neq\infty$ assume that the $\ell_i$ are
different from $p$.
Let $n=\ell_1\cdot \cdots \cdot \ell_r$. Then
$$\frac{\vert S^{(q)}(n)\vert}{\vert\GSp^{(q)}_{2g}(\Zz/n\Zz)\vert}=\prod_{j=1}^r\frac{\vert S^{(q)}(\ell_j)\vert}
{\vert \GSp_{2g}^{(q)}(\F_{\ell_j})\vert}$$ where $S^{(q)}(n)\subset
\GSp_{2g}(\Zz/n\Zz)$ is the set of matrices that belong to
$S^{(q)}(\ell_j)$ modulo $\ell_j$, for all $j=1, \dots, r$.
\end{itemize}
\end{thm}

First we will prove part (1) of Theorem \ref{Proposicion}. We will need four
lemmata.

On the one hand, the cardinality of $S_{\lambda}(\ell)_0$ is very
easy to compute.

\begin{lemm}\label{Greater} It holds that
$$\vert  S_{\lambda}(\ell)_0\vert =
\ell^{2g-2}(\ell-1)\beta(\ell, g-1)\vert\Sp_{2g-4}(\F_{\ell})\vert.$$
\end{lemm}

Moreover we can compute the cardinality of $S_{\lambda}(\ell)$ in
terms of $\vert S_{\lambda}(\ell)_0\vert$.

\begin{lemm}\label{Multiplicatividad} $\vert S_{\lambda}(\ell)\vert=(\ell^{2g-2}(\ell-1)+1)\vert S_{\lambda}(\ell)_0\vert$.
\end{lemm}

{\em Proof.} Let $A\in S_{\lambda}(\ell)_0$. First of all we will
see that the vectors fixed by $A$ are those in the cyclic subgroup
generated by $e_1$. Since the matrix $A$ clearly fixes the vectors
in the cyclic subgroup generated by $e_1$, it suffices to show that
any vector fixed by $A$ must belong to this subgroup.

Consider the system of equations $A(x_1, \dots, x_{2g})^t=(x_1,
\dots, x_{2g})^t$. Assume first that we have a solution with
$x_2=0$. Then the last $2g-2$ equations boil down to
$$B(x_3, \dots, x_{2g})^t=(x_3, \dots, x_{2g})^t.$$ But since $B$ does not have the eigenvalue $1$, this equations are
not simultaneously satisfied by a nonzero tuple, hence $(x_1, \dots,
x_{2g})^t$ belongs to the cyclic group generated by $e_1$.

Assume now that we have a solution $(x_1, \dots, x_g)^t$ with
$x_2\not=0$. Since $1$ is not an eigenvalue of $B$, the matrix
$\Id-B$ is invertible, and we can write the last $2g-2$ equations as

$$(x_3/x_2,
\dots,  x_{2g}/x_2)^t=(\Id-B)^{-1}(d_1,  \dots,  d_{2g-2})^t.$$ On
the other hand, the first equation reads
$$d=-(b_1, \dots, b_{2g-2})(x_3/x_2, \dots, x_{2g}/x_2)^t.$$

Hence

$$d=-(b_1, \dots, b_{2g-2})(\Id-B)^{-1}(d_3,
\dots, d_{2g})^t.$$ But we have precisely asked that $d$ does not
satisfy such an equation, cf. Definition \ref{Definicion}.

Now one can check that, if we have $A, \widetilde{A}\in
S_{\lambda}(\ell)_0$ and elements $\alpha_3, \dots, \alpha_{2g}$,
$\beta, \widetilde{\alpha}_{3}, \dots, \widetilde{\alpha}_{2g}$,
$\widetilde{\beta}$ in $\F_{\ell}$ such that
$T_{u_{\alpha}}[\beta]^{-1} \cdot A\cdot
T_{u_{\alpha}}[\beta]=T_{u_{\widetilde{\alpha}}}[\widetilde{\beta}]^{-1}
\cdot \widetilde{A}\cdot
T_{u_{\widetilde{\alpha}}}[\widetilde{\beta}]$, then either
$\beta=\widetilde{\beta}=0$ and $A=\widetilde{A}$ or else
$\alpha_k=\widetilde{\alpha}_k$ for $k=3, \dots, 2g$,
$\beta=\widetilde{\beta}$ and $A=\widetilde{A}$. Namely, one notices
that since
$\widetilde{A}=T_{u_{\widetilde{\alpha}}}[\widetilde{\beta}]
T_{u_{\alpha}}[\beta]^{-1} \cdot A\cdot
T_{u_{\alpha}}[\beta]T_{u_{\widetilde{\alpha}}}[\widetilde{\beta}]^{-1}$
fixes $e_1$, then $A$ fixes
$T_{u_{\alpha}}[\beta]T_{u_{\widetilde{\alpha}}}[\widetilde{\beta}]^{-1}e_1$.
But $A$ only fixes the elements of the cyclic group generated by
$e_1$; hence,
$T_{u_{\alpha}}[\beta]T_{u_{\widetilde{\alpha}}}[\widetilde{\beta}]^{-1}e_1$
must be in the cyclic group generated by $e_1$. Now  computing
$T_{u_{\alpha}}[\beta]T_{u_{\widetilde{\alpha}}}[\widetilde{\beta}]^{-1}e_1$
one can conclude easily.

Therefore each element of $S_{\lambda}(\ell)_0$ gives rise to a
subset of $S_{\lambda}(\ell)$  by conjugation by the matrices
$T_{u_{\alpha}}[\beta]$, where $\alpha$ runs through the tuples
$(\alpha_3, \dots, \alpha_{2g})\in \F_{\ell}^{2g-2}$ and $\beta\in
\F_{\ell}$, and $S_{\lambda}(\ell)$ is the disjoint union of these
subsets. Furthermore, each of these sets has cardinality
$\ell^{2g-2}(\ell-1) +1$.

\hfill $\Box$

To prove the first part of Theorem \ref{Proposicion}, we only need one more
lemma, which is an easy consequence of the Chinese Remainder
Theorem.

\begin{lemm}\label{Disjointness} \begin{itemize}
\item [(1)] Let $q$ be a power of a prime number $p$, and let $n$ be a
squarefree natural number such that $p\nmid n$. The cardinality of
$\GSp^{(q)}_{2g}(\Zz/n\Zz)$ equals $
\ord_{n}(q)\cdot\prod_{\ell\vert n}\vert \Sp_{2g}(\F_{\ell})\vert$.
\item [(2)] Let $q=\infty$, and let $n$ be a squarefree natural number.
The cardinality of $\GSp^{(q)}_{2g}(\Zz/n\Zz)$ equals $
\prod_{\ell\vert n}(\ell-1)\vert \Sp_{2g}(\F_{\ell})\vert$.
\end{itemize}
\end{lemm}

{\em Proof of Theorem \ref{Proposicion}(1)}

Let $q$ be a power of a prime $p$ or $q=\infty$, and let $\ell$ be a
prime. In the first case, let us also assume $\ell\not=p$. Applying
the identity \eqref{Induccion} in the proof of Lemma \ref{Independent} twice to the
cardinality of $\GSp_{2g}^{(q)}(\F_{\ell})$ and Lemmas \ref{Multiplicatividad},
\ref{Greater} and \ref{Disjointness}, we obtain
$$\frac{\vert S^{(q)}(\ell)\vert}{\vert \GSp_{2g}^{(q)}(\F_{\ell})\vert}=\frac{(\ell^{2g-2}(\ell-1)+1)\ell^{2g-2}(\ell-1)\beta(\ell, g-1)\vert
\Sp_{2g-4}(\F_{\ell})\vert}{(\ell^{2g}-1)\ell^{2g-1}
(\ell^{2g-2}-1)\ell^{2g-3}\vert \Sp_{2g-4}(\F_{\ell})\vert}\sim
\frac{1}{\ell},$$ and the sum $\sum_{\ell\not=p\text{ prime
}}\frac{1}{\ell}$ diverges.\hfill $\Box$

To prove the second part of Theorem \ref{Proposicion} we need one
auxiliary lemma. For the rest of the section, $q$ will be a power of
a prime $p$.

For each squarefree $n$ not divisible by $p$ and each $i=1, \dots,
\ord_n(q)$, define $S_{q^i}(n):=\{A\in S^{(q)}(n): \varepsilon(A)=q^i
\text{ modulo }n\}$.

\begin{lemm}\label{Previous}  Let $q$ be a power of a prime number $p$.
Let $\ell_1, \dots, \ell_r$ be distinct primes which are different from $p$, and consider
$n=\ell_1\cdot\cdots\cdot \ell_r$. Let $i\in \{1, \dots, \ord_n(q)\}$. Then
there is a bijection
$$S_{q^i}(n)\simeq S_{q^i}(\ell_1)\times \cdots \times S_{q^i}(\ell_r).$$
\end{lemm}

{\em Proof.}
Consider the canonical projection
$$\begin{aligned}\pi:S_{q^i}(n)&\rightarrow S_{q^i}(\ell_1)\times \cdots \times S_{q^i}(\ell_r)\\
A& \mapsto (A\mod{\ell_1}, \dots, A\mod{\ell_r}).\end{aligned}$$

This is clearly an injective map. Now we want to prove surjectivity.
For each $j$, take some matrix $B_j\in S_{q^i}(\ell_j)$.

By the Chinese Remainder Theorem, there exists  $A\in
\GSp_{2g}(\Zz/n\Zz)$ such that $A$ projects onto $B_j$ for each $j$.
Note that in particular $A\in S^{(q)}(n)$. Since $\varepsilon(A)$ is
congruent to $\varepsilon(B_j)=q^i$ modulo $\ell_j$ for all $j$, we
get that $\varepsilon(A)=q^i$ modulo $n$. Therefore $A\in
S_{q^i}(n)$. \hfill $\Box$

{\em Proof of Theorem \ref{Proposicion}(2)}\\
{\em Case $q\neq\infty$:} On the one hand, since the cardinality of $\vert S_{q^i}(\ell)\vert$
does not depend on $i$, we obtain

$$\prod_{\ell\vert n}\frac{\vert S^{(q)}(\ell)\vert}{\vert \GSp_{2g}^{(q)}(\F_{\ell})\vert}=
\prod_{\ell\vert n} \frac{\ord_{\ell}(q) \vert
S_q(\ell)\vert}{\ord_{\ell}(q)\vert \Sp_{2g}(\F_{\ell})\vert}=
\prod_{\ell\vert n} \frac{\vert S_q(\ell)\vert}{\vert
\Sp_{2g}(\F_{\ell})\vert}.$$

On the other hand, taking into account again that $\vert
S_{q^i}(\ell)\vert$ is independent of $i$, Lemma \ref{Disjointness},
and that $\vert S_{q^i}(n)\vert=\prod_{\ell\vert n}\vert
S_{q^i}(\ell)\vert$ by Lemma \ref{Previous}, we get

$$\frac{\vert S^{(q)}(n)\vert}{\vert \GSp_{2g}^{(q)}(\Zz/n\Zz)\vert}=
\frac{\sum_{i=1}^{\ord_n(q)}\vert
S_{q^i}(n)\vert}{\ord_n(q)\prod_{\ell\vert n}\vert
\Sp_{2g}(\F_{\ell})
\vert}=\frac{\sum_{i=1}^{\ord_n(q)}\left(\prod_{\ell\vert n}\vert
S_{q^i}(\ell)\vert\right)}{\ord_n(q)\prod_{\ell\vert n}\vert
\Sp_{2g}(\F_{\ell}) \vert}=$$
$$=\frac{\sum_{i=1}^{\ord_n(q)}\left(\prod_{\ell\vert n}\vert
S_q(\ell)\vert\right)}{\ord_n(q)\prod_{\ell\vert n}\vert
\Sp_{2g}(\F_{\ell}) \vert}=\frac{\ord_{n}(q)\left(\prod_{\ell\vert
n}\vert S_q(\ell)\vert\right)}{\ord_n(q)\prod_{\ell\vert n}\vert
\Sp_{2g}(\F_{\ell}) \vert}=$$
$$=\prod_{\ell\vert n} \frac{\vert S_q(\ell)\vert}{\vert \Sp_{2g}(\F_{\ell})\vert}.$$

{\em Case $q=\infty$:} By the Chinese Remainder Theorem, there is a canonical isomorphism
$$c: \GSp^{(\infty)}_{2g}(\Zz/n) \cong \prod_{i=1}^r \GSp^{(\infty)}_{2g}(\Zz/\ell_i)$$
and $$S^{(\infty)}(n)=c^{-1}(S^{(\infty)}(\ell_1)\times\cdots\times S^{(\infty)}(\ell_r))$$
by the definition of $S^{(\infty)}(n)$.
It follows that
$$\frac{|S^{(\infty)}(n)|}{|\GSp^{(\infty)}_{2g}(\Zz/n)|}
=\prod_{i=1}^r \frac{|S^{(\infty)}(\ell_i)|}{|\GSp^{(\infty)}_{2g}(\Zz/\ell_i)|}$$
as desired.\hfill $\Box$

\begin{rema} In the definition of the set $S_{q^i}(\ell)_0$ (cf. Definition \ref{Definicion}), we
choose
a subset $\calB_{q^i}$ of matrices in $\GSp_{2g-2}(\F_{\ell})[q^i]$
without the eigenvalue $1$, which is large enough to ensure that
part (1) of Theorem \ref{Proposicion} holds. For a concrete
value of $g$, one can choose such set more explicitly. For instance,
when $g=2$, instead of $\calB_{q^i}$ one can consider the set
\begin{multline*}\calB'_{q^i}:=\{\begin{pmatrix}b_{1, 1} & b_{1, 2}\\ b_{2, 1} & b_{2, 2}\end{pmatrix}:
b_{1, 1}\in \F_{\ell},\ b_{2, 2}\in \F_{\ell}\setminus\{1-b_{1, 1} +
q^i\},\\ b_{1, 2}\in \F_{\ell}^\times,\ b_{2, 1}=b_{1, 2}^{-1}(b_{1,
1}b_{2, 2} - q^i)\}\end{multline*} of $\ell(\ell-1)^2$ matrices, which can also be used to prove the second part of Theorem \ref{Proposicion} in the case of the group $\GSp_4(\F_{\ell})$.

\end{rema}

\section{Proof of the Conjecture of Geyer and Jarden, part a)}

\begin{thm} \label{GJmain2} Let $(A, \lambda)$ be a
polarized abelian variety over a finitely generated field $K$.
Assume that $A/K$ has big monodromy. Then for almost all
$\sigma\in G_K$ there are infinitely many prime
numbers $\ell$ such that $A(K\sep(\sigma))[\ell]\neq 0$.
\end{thm}

{\em Proof.} Let $p:= \chara(K)$. Let $G=G_K$ and $g:=\dim(A)$.
We fix once and for all for every prime number
$\ell\neq p$ a symplectic basis of $T_\ell A$. This defines an isometry
of symplectic spaces $(A[n], e_n^\lambda)\cong ((\Zz/n)^{2g},
e_n^{\can})$, where $e_n^{\can}$ denotes the standard symplectic
pairing on $(\Zz/n)^{2g}$, for every $n\in \Nn$ which is not
divisible by $p$. We get an isomorphism
$\GSp(A[n], e_n^\lambda)\cong \GSp_{2g}(\Zz/n)$ for
every such $n$, and we consider the representations
$$\rho_n: G_K\to \GSp_{2g}(\Zz/n)$$
attached to $A/K$ after these choices. If $m$ is a divisor of $n$, then
we denote by $r_{n, m}: \GSp_{2g}(\Zz/n)\to
\GSp_{2g}(\Zz/m)$ the corresponding canonical map, such that
$r_{n, m}\circ \rho_n=\rho_m$.

Let $q:=q(K)$ be the cardinality of the algebraic closure
of the prime field of $K$ in $K$. Thus $q=\infty$ if $p=0$ and $q$ is a power of $p$ otherwise.
As $A$ has big monodromy, we find by Proposition \ref{fullresult} an
integer $c$ (divisible by $p$, if $p\neq 0$) such that $\im(\rho_n)=\GSp^{(q)}_{2g}(\Zz/n),$
for every $n$ coprime to $c$.

For every prime number $\ell> c,$ we define
$$X_\ell:=\{\sigma\in G_K\,|\, A(K\sep(\sigma))[\ell]\neq 0\}.$$
It is enough to prove that $\bigcap_{n>c}\bigcup_{\ell\ge n \
\mathrm{prime}} X_\ell$ has measure $1$. Let
$S^{(q)}(n)\subset\GSp_{2g}^{(q)}(\Zz/n)$ be the special sets of
symplectic matrices defined in Section 4. By remark \ref{eveins}
$\rho_\ell^{-1}(S^{(q)}(\ell))\subset X_\ell$ for every prime number
$\ell>c$. It is thus enough to prove that
$\bigcap_{n>c}\bigcup_{\ell\ge n \ \mathrm{prime}}
\rho_\ell^{-1}(S^{(q)}(\ell))$ has measure $1$. Clearly
$\mu_G(\rho_n^{-1}(S^{(q)}(n)))=\frac{|S^{(q)}(n)|}{|\GSp_{2g}^{(q)}(\Zz/n)|}$
for all integers $n$ coprime to $c$. Hence part (1) of Theorem
\ref{Proposicion} implies that $\sum_{\ell>c\ \mathrm{prime}}
\mu_G(\rho_\ell^{-1}(S^{(q)}(\ell)))=\infty$.

Furthermore, if $\ell_1,\cdots,\ell_r>c$ are different prime numbers and $n=\ell_1\cdots \ell_r$, then
$$\bigcap_{i=1}^r \rho_{\ell_i}^{-1}(S^{(q)}(\ell_i))=\rho_n^{-1}(S^{(q)}(n))$$
and part (2) of Theorem \ref{Proposicion} implies
$$\mu_G( \bigcap_{i=1}^r \rho_{\ell_i}^{-1}(S^{(q)}(\ell_i)))=
\prod_{i=1}^r \mu_G(\rho_{\ell_i}^{-1}(S^{(q)}(\ell_i))).$$
Hence $(\rho_{\ell}^{-1}(S^{(q)}(\ell)))_{\ell>c}$ is a $\mu_G$-independent sequence of subsets of $G$.
 Now Lemma \ref{borelcantelli} implies that $\bigcap_{n>c}\bigcup_{\ell\ge n\ \mathrm{prime}} \rho_\ell^{-1}(S^{(q)}(\ell))$  has
measure $1$, as desired.\hfill $\Box$
\bigskip\bigskip

\noindent
{\bf Acknowledgements.}
S. A. is a research fellow of the Alexander von Humboldt Foundation.
S. A. was partially supported by the Ministerio de Educaci\'on y Ciencia grant MTM2009-07024. S. A. wants to thank the Hausdorff Research Institute for Mathematics in Bonn, the Centre de Recerca Matem\`atica in Bellaterra and the Mathematics Department of Adam Mickiewicz University in Pozna\'n for their support and hospitality while she worked on this project.
W.G. and S.P. were partially supported by the Deutsche Forschungsgemeinschaft
research grant GR 998/5-1.
W.G. was partially supported by the Alexander von Humboldt Research Fellowship
and an MNiSzW grant. W.G. thanks Centre Recerca Matem\`atica in Bellaterra and the Max Planck
Institut f{\" u}r Mathematik in Bonn for support and hospitality during visits in 2010, when he worked on this
project. S.P. gratefully acknowledges the hospitality of Mathematics Department of Adam Mickiewicz  University in Pozna{\' n} and of the Minkowski center at Tel Aviv University during several research visits.


\end{document}